\documentclass[reqno,12pt]{article} \usepackage{fontenc}
\usepackage{amsfonts}
\usepackage[utf8]{inputenc}
\usepackage[english]{babel}
\usepackage{amsmath}
\usepackage{amssymb}

\newtheorem{theorem}{Theorem}
\newtheorem{corollary}{Corollary}[theorem]

\newtheorem{example}{Example}[theorem]
\newtheorem{definition}{Definition}[theorem]
\begin{document}
{\small {~~\\ \vspace{-3.5cm} \begin{flushright}\jobname \\ \today
\end{flushright} \vspace{1.5cm} }}

\begin{center}
{\bf \Large
Restricted Hom-Lie Superalgebras}\\~~\\
{Shadi Shaqaqha}\\~\\
Yarmouk University, Irbid, Jordan\\
shadi.s@yu.edu.jo
\end{center}

%
%
%

\begin{abstract}
The aim of this paper is to introduce the notion of restricted Hom-Lie superalgebras. This class of algebras is a generalization of both restricted Hom-Lie algebras and restricted Lie superalgebras. In this paper, we present a way to obtain restricted Hom-Lie superalgebras from the classical restricted Lie superalgebras along with algebra endomorphisms. Homomorphisms relations between restricted Hom-Lie superalgebras are defined and studied. Also, we obtain some properties of $p$-maps and restrictable Hom-Lie superalgebras.
\end{abstract}

\vspace{.4cm} \noindent {\bf Keywords:} Hom-Lie superalgebra; Restricted Hom-Lie superalgebra; $p$-map; Restrictable Hom-Lie superalgebra; Morphism.

\section{INTRODUCTION}

Hom structures including Hom-algebras, Hom-Lie algebras, Hom-Lie superalgebras, Hom-Lie color algebras, Hom-coalgebras,
Hom-modules, and Hom-Hopf modules have been widely investigated during the last years. The motivations to study Hom-Lie structures are related to physics and to deformations of Lie algebras, especially Lie algebras of
vector fields.  The Hom-Lie algebras were firstly studied by Hartwig, Larsson, and Silvestrov in \cite{Hartwig}. Later, Hom-Lie superalgebras are introduced by Ammar and Makhlouf in \cite{Makhlouf}. Also, Hom-Lie color algebras, which are the natural generalizations of Hom-Lie algebras and Hom-Lie superalgebras, are studied by Yuan (\cite{Yuan}). However, the notion of restricted Hom-Lie algebras was introduced by Guan and Chen (\cite{Guan}).\\
The Lie superalgebras were originally introduced by Kac (\cite{Kac}). Later, when $\mathrm{char}F=p > 0$, Lie superalgebras were studied, and also the notion of restricted Lie superalgebras was introduced (\cite{Kochetkov, Petrogradsky3}). In the present work, we introduce the notion of restricted Hom-Lie superalgebras. \\
The article is organized as follows. In Section 2, we review some definitions, notations, and results in \cite{Makhlouf, Bahturin 3}. In Section 3, we introduce the definition of restricted Hom-Lie superalgebra. We provide some properties and their relationships with restricted Hom-Lie superalgebras. As a result, we show that it is enough to know how a $p$-map behaves for inputs from any basis of the domain. In Section 4, We study the direct sum of (restricted) Hom-Lie superalgebras. In Section 5, we study the homomorphisms between restricted Hom-Lie superalgebras. In particular, we show how arbitrary restricted Lie superalgebras deform into restricted Hom-Lie superalgebras via algebra endomorphisms. Also we discuss the images as well as
preimages of restricted Hom-Lie subsuperalgebras under homomorphisms. The restrictable Hom-Lie superalgebras are defined and studied in Section 6.

\section{Preliminaries}

 Let $F$ be the ground field of characteristic $\neq 2, 3$. A linear superspace $V$ over $F$ is merely a $\mathbb{Z}_2$-graded linear space with a direct sum $V=V_0\oplus V_1$. The elements of $V_j$, $j\in \{0, 1\}$, are said to be homogeneous of parity $j$. The parity of
a homogeneous element $x$ is denoted by $|x|$. Suppose that $V=V_0\oplus V_1$ and $V'=V_0'\oplus V_1'$ are two linear superspaces. A linear map $\alpha:V\rightarrow V'$ is an even if $\alpha(V_j)\subseteq V_j'$ for $j\in\{0,1\}$.

\begin{definition}[{{\cite{Makhlouf}}}]\label{CB}
A Hom-associative superalgebra is a triple $(A, \mu, \alpha)$ where $A$ is a linear superspace, $\mu:A\times A\rightarrow A$ is an even bilinear map, and $\alpha:A\rightarrow A$ is an even linear map such that
$$\mu(\alpha(x), \mu(y,z))=\mu(\mu(x, y), \alpha(z)).$$
\end{definition}
A Hom-associative superalgebra $(A, \mu, \alpha)$ is called multiplicative if $\alpha(\mu(x, y))=\mu(\alpha(x), \alpha(y))$ for all $x, y\in A$.
\begin{definition}[{{\cite{Makhlouf}}}] \label{CBA}
A Hom-Lie superalgebra is a triple $(L, [~,~], \alpha)$ where $L$ is a linear superspace, $[~,~]:L\times L\rightarrow L$ is an even bilinear map, and $\alpha:L\rightarrow L$ is an even linear map such that the following identities satisfied for any homogeneous $x, y, z\in L$.
\begin{itemize}
\item[(i)] Super skew-symmetry:
$$[x, y]= -(-1)^{|x||y|}[y,x].$$
\item [(ii)] Hom-superJacobi identity:
$$(-1)^{|x||z|}[\alpha(x), [y,z]]+(-1)^{|z||y|}[\alpha(z), [x, y]]+(-1)^{|y||x|}[\alpha(y), [z,x]]=0.$$
\end{itemize}
\end{definition}
\begin{example}[{{\cite{Makhlouf}}}]\label{C}
Any Hom-associative superalgebra $A$, with the bracket
$$[x, y]=\mu(x, y)-(-1)^{|x||y|}\mu(y, x)$$
for any nonzero homogeneous $x, y\in A$, is a Hom-Lie superalgebra, which will be denoted by $A^{(-)}$.
\end{example}
It is clear that Lie superalgebras are examples of Hom-Lie superalgebras by setting $\alpha=\mathrm{id}_L$. A Hom-Lie superalgebra is called a multiplicative Hom-Lie superalgebra if $\alpha$ is an even homomorphism (that is, $\alpha([x, y])=[\alpha(x), \alpha(y)]$ for all $x, y\in L$). A subspace $H\subseteq L$ is called a Hom-Lie subsuperalgebra if $\alpha(H)\subseteq H$ and $H$ is closed under the bracket operation (that is, $[x, y]\in L$ for all $x, y\in L$).\\
Let $L=L_0\oplus L_1$ be a Lie superalgebra. For a homogenous element $a\in L$, we consider a mapping $\mathrm{ad}a:L\rightarrow L;~ b\mapsto [a, b]$.

\begin{definition}[{{\cite{Bahturin 3}}}] \label{A1}
A Lie superalgebra $L=L_0\oplus L_1$ is called {\em restricted} (or {\em $p$-superalgebra}) if there is a map
$$^{[p]}~:L_0\rightarrow L;~ x\mapsto x^{[p]},$$
satisfying
\begin{itemize}
\item[(i)] $\left(\mathrm{ad}x\right)^p= \mathrm{ad}\left(x^{[p]}\right)$ for all $x\in L_0$,
\item[(ii)] $\left(k x\right)^{[p]}= k^px^{[p]}$ for all $k \in F, x\in L_0$,
\item[(iii)] $\left(x+y\right)^{[p]}=x^{[p]} + y^{[p]}+\sum_{i=1}^{p-1}s_i(x, y)$ for all $x, y\in L_0$ where $is_i(x, y)$ is the coefficient of $\lambda^{i-1}$ in $\mathrm{ad}\left(\lambda x+ y\right)^{p-1}(x)$.
\end{itemize}
\end{definition}
\section{Restricted Hom-Lie superalgebra}
Let $(L, [~,~], \alpha)$ be a multiplicative Hom-Lie superalgebra. For a homogeneous element $a\in L_0$ with $\alpha(a)= a$, we consider a map
$$\mathrm{ad}_{\alpha}a:L\rightarrow L;~b\mapsto [a, \alpha(b)].$$
Put $L^0= \left\{x\in L_0~|~\alpha(x)= x\right\}$. We can easily prove that $L^0$ is a Hom-Lie subsuperalgebra of $L$.
\begin{definition}\label{A2}
A multiplicative Hom-Lie superalgebra $(L,[~,~], \alpha)$ is called restricted if there is a map (called a $p$-map)
$$[p]~:L^0\rightarrow L^0;~x\mapsto x^{[p]},$$
satisfying
\begin{itemize}
\item[(i)] $\left(\mathrm{ad}_{\alpha}x\right)^p=\mathrm{ad}_{\alpha}\left(x^{[p]}\right)$ for all $x\in L^0$,
\item[(ii)] $\left(k x\right)^{[p]}= k^px^{[p]}$ for all $k \in F, x\in L^0$,
\item[(iii)] $\left(x+y\right)^{[p]}=x^{[p]} + y^{[p]}+\sum_{i=1}^{p-1}s_i(x, y)$ for all $x, y\in L^0$ where $is_i(x, y)$ is the coefficient of $\lambda^{i-1}$ in $\mathrm{ad}_{\alpha}\left(\lambda x+ y\right)^{p-1}(x)$.
\end{itemize}
\end{definition}
Let $(L, [~,~], \alpha, [p])$ be a restricted Hom Lie-superalgebra over a field $F$. A Hom-Lie subsuperalgebra $H=H_0\oplus H_1$ of $L$ is called a $p$-subalgebra if $x^{[p]}\in H^0$ $\forall x\in H^0$.
\begin{definition}\label{A3}
Let $(L, [~,~], \alpha)$ be a Hom Lie-superalgebra, and let $S$ be a subset of $L$.
\begin{itemize}
\item[(i)] A map $f:L\rightarrow L$ is called a $p$-semilinear map if $f(kx + y)= k^pf(x)+ f(y)$ $\forall x, y\in L$ and $\forall k\in F$.
\item[(ii)] The $\alpha$-centralizer of $S$ in $L$, denoted by $C_L(S)$, is defined to be
$$C_L(S)=\left\{x\in L~|~[x, \alpha(y)]=0~\forall y\in S\right\}.$$
In particular, if $S=L$, it is called the $\alpha$-center of $L$, and it is denoted by $C(L)$.
\end{itemize}
\end{definition}
The following theorem was proved by Baoling Guan and Liangyun Chen in \cite{Guan} in the setting of Hom-Lie algebras.
\begin{theorem}\label{A4}
Let $H$ be a Hom-Lie subsuperalgebra of a multiplicative restricted Hom-Lie superalgebra $(L, [~,~], \alpha, [p])$ and $[p]_1:H^0\rightarrow H^0$ is a map. Then the following statements are equivalent:
\begin{itemize}
\item[(i)] $[p]_1$ is a $p$-map on $H^0$,
\item[(ii)] there is a $p$-semilinear map $f: H^0\rightarrow C_L(H)$ such that $[p]_1=[p]+f$.
\end{itemize}
\end{theorem}
{\it Proof.~}
Suppose that $[p]_1$ is a $p$-map on $H^0$. Consider
$$f: H^0\rightarrow L;~x\mapsto x^{[p]_1}- x^{[p]}.$$
Since $\mathrm{ad}_\alpha f(x)(y)=[f(x), \alpha(y)]=[x^{[p]_1}, \alpha(y)]-[x^{[p]}, \alpha(y)]= \left(\mathrm{ad}_{\alpha}x\right)^p-\left(\mathrm{ad}_{\alpha}x\right)^p=0$ for all $x\in H^0, y\in L$, $f$ actually maps $H^0$ into $C_L(H)$. Now, for $x, y\in H^0$ and $k\in F$, we obtain
\begin{eqnarray*}
f(kx+y)&=&(kx+y)^{[p]_1}-(kx+y)^{[p]}\\
&=& k^px^{[p]_1} +y^{[p]_1}+\sum_{i=1}^{p-1}s_i(kx, y)-k^px^{[p]}-y^{[p]}- \sum_{i=1}^{p-1}s_i(kx, y)\\
&=&k^pf(x)+f(y).
\end{eqnarray*}
This shows that $f$ is a $p$-semilinear map. Conversely, assume there exists a $p$-semilinear map $f: H^0\rightarrow C_L(H)$ with $[p]_1=[p]+f$. We check the three conditions given in Definition \ref{A2}. For $x\in H^0$, $y\in H$, we have
\begin{eqnarray*}
\mathrm{ad}_\alpha x^{[p]_1}(y)&=&\mathrm{ad}_\alpha(x^{[p]}+f(x))(y)\\
                              &=&\mathrm{ad}_\alpha x^{[p]}(y)+\mathrm{ad}_\alpha f(x) (y)\\
                              &=&\mathrm{ad}_\alpha x^{[p]}(y)~~(\mathrm{since}~f(x)\in C_L(H))\\
                              &=&(\mathrm{ad}_\alpha x)^p (y).
\end{eqnarray*}
For $x, y\in H^0$, we have
\begin{eqnarray*}
(x+y)^{[p]_1}&=&(x+y)^{[p]}+f(x+y)\\
&=& x^{[p]}+y^{[p]}+\sum_{i=1}^{p-1}s_i(x, y)+f(x)+f(y)\\
&=&x^{[p]_1}+y^{[p]_1}+\sum_{i=1}^{p-1}s_i(x, y).
\end{eqnarray*}
and, for $k\in F$, we get
\begin{eqnarray*}
(kx)^{[p]_1}&=&(kx)^{[p]}+f(kx)\\
&=&k^px^{[p]}+k^pf(x)\\
&=&k^p(x^{[p]}+f(x))\\
&=&k^px^{[p]_1}.
\end{eqnarray*}
\hfill $\Box$
\begin{corollary}\label{A5}
Let $(L, [~,~], \alpha)$ be a multiplicative Hom-Lie superalgebra.
\begin{itemize}
\item[(i)] If $C(L)=0$, then $L$ admits at most one $p$-map.
\item[(ii)] If two $p$-maps coincide on a basis of $L^0$, then they are equal.
\item[(iii)] If $[p]: L^0\rightarrow L^0$ is a $p$-map, then there exists a $p$-mapping $[p]'$ of $L$ such that $x^{[p]'}=0$ $\forall x\in C(L^0)$.
\end{itemize}
\end{corollary}
{\it Proof.~}
\begin{itemize}
\item[(i)] Suppose $[p]_1$ and $[p]_2$ are $p$-maps of $L$. By Theorem \ref{A4}, there exists a $p$-semilinear map $f$ from $L^0$ into $C(L)$ such that $[p]_2=[p]_1+f$, and since $C(L)=0$, we have $[p]_1=[p]_2$.
\item[(ii)] Let $[p]_1$ and $[p]_2$ be two $p$-maps coincide on a basis $B$ of $L$. According to Theorem \ref{A4}, there exists a $p$-semilinear map $f:L^0\rightarrow C(L)$ such that $f=[p]_2-[p]_1$, and so $f(b)=0$ $\forall b\in B$. As $f$ is a $p$-semilinear, we have $f(x)=0$ $\forall x\in L^0$. Thus, $[p]_1=[p]_2$.
\item[(iii)] $[p]|_{C(L^0)}$ is obviously a $p$-map. For $x, y\in C(L^0)$ and $k\in F$, we have
    \begin{eqnarray*}
    (kx+y)^{[p]}&=&(kx)^{[p]}+y^{[p]}+ \sum_{i=1}^{p-1}s_i(kx, y)\\
                &=&k^px^{[p]}+y^{[p]} ~~~~~~(\mathrm{since ~C(L^0)~ is~abelian}).
    \end{eqnarray*}
Thus it is a $p$-semilinear map. Extend it to a $p$-semilinear map $f:L^0\rightarrow C(L^0)$. Then $[p]':=[p]-f$ is a $p$-map with $x^{[p]'}=0$ $\forall x\in C(L^0)$.
\end{itemize}
\hfill $\Box$

\section{Direct sum of Hom-Lie superalgebras}

Let $V$ and $V'$ be two superspace. Then the homogeneous elements of the superspace $V\oplus V'$ have the form $(x,y)$, where $x$ and $y$ are homogeneous elements of the same degrees in $V$ and $V'$, respectively. In this case, $|(x, y)|=|x|=|y|$. We have the following result.
\begin{theorem}\label{A8}
Given two Hom-Lie superalgebras $(L, [~,~]_L, \alpha)$ and $(\Gamma, [~,~]_\Gamma, \beta)$, then $(L\oplus \Gamma, [~, ~]_{L\oplus \Gamma}, \gamma)$ has a multiplicative Hom-Lie superalgebra structure, where the bilinear map $[~, ~]_{L\oplus \Gamma}$ is given by ($(x_1, y_1)$ and $(x_2,y_2)$ are homogeneous elements)
$$[(x_1, y_1), (x_2, y_2)]_{L\oplus \Gamma}=([x_1, x_2]_L,[y_1, y_2]_\Gamma),$$
and the linear map $\gamma:L\oplus \Gamma\rightarrow L\oplus \Gamma$ is given by
$$\gamma(x_1,y_1)=(\alpha(x_1),\beta(y_1)).$$
Moreover, $(L\oplus \Gamma)^0=L^0\oplus \Gamma^0$.
\end{theorem}
{\it Proof.~}
Suppose that $(x_1, y_1)$ and $(x_2, y_2)$ are two homogeneous elements in $L\oplus \Gamma$. Then
\begin{eqnarray*}
[(x_1, y_1), (x_2, y_2)]_{L\oplus \Gamma}&=& ([x_1, x_2]_L, [y_1, y_2]_{\Gamma})\\
                                         &=&(-(-1)^{|x_1||x_2|}[x_2, x_1]_L, -(-1)^{|y_1||y_2|}[y_2, y_1]_{\Gamma})\\
                                         &=&-(-1)^{|x_1||x_2|}([x_2, x_1]_L, [y_2, y_1]_{\Gamma})~~~(|x_1|=|y_1|~\mathrm{and}~|x_2|=|y_2|)\\
                                         &=&-(-1)^{|(x_1, y_1)||(x_2, y_2)|}[(x_2, y_2), (x_1, y_1)]_{L\oplus \Gamma}.
\end{eqnarray*}
For all homogeneous elements $(x_1, y_1), (x_2, y_2)$, and $(x_3, y_3)$, we have,
\begin{eqnarray*}
&&(-1)^{|(x_1, y_1)||(x_3, y_3)|}[\gamma(x_1, y_1), [(x_2, y_2), (x_3, y_3)]_{L\oplus \Gamma}]_{L\oplus \Gamma}\\
&+&(-1)^{|(x_3, y_3)||(x_2, y_2)|}[\gamma(x_3, y_3), [(x_1, y_1), (x_2, y_2)]_{L\oplus \Gamma}]_{L\oplus \Gamma}\\
&+&(-1)^{|(x_2, y_2)||(x_1, y_1)|}[\gamma(x_2, y_2), [(x_3, y_3), (x_1, y_1)]_{L\oplus \Gamma}]_{L\oplus \Gamma}\\
&=&(-1)^{|x_1||x_3|}[(\alpha(x_1), \beta(y_1)), ([x_2, x_3]_L, [y_2, y_3]_{\Gamma})]_{L\oplus \Gamma}\\
&+&(-1)^{|x_3||x_2|}[(\alpha(x_3), \beta(y_3)), ([x_1, x_2]_L, [y_1, y_2]_{\Gamma})]_{L\oplus \Gamma}\\
&+&(-1)^{|x_2||x_1|}[(\alpha(x_2), \beta(y_2)), ([x_3, x_1]_L, [y_3, y_1]_{\Gamma})]_{L\oplus \Gamma}\\
&=&((-1)^{|x_1||x_3|}[\alpha(x_1), [x_2, x_3]_L]_L, (-1)^{|y_1||y_3|}[\beta(y_1), [y_2, y_3]_\Gamma]_\Gamma)\\
&+&((-1)^{|x_3||x_2|}[\alpha(x_3), [x_1, x_2]_L]_L, (-1)^{|y_3||y_2|}[\beta(y_3), [y_1, y_2]_\Gamma]_\Gamma)\\
&+&((-1)^{|x_2||x_1|}[\alpha(x_2), [x_3, x_1]_L]_L, (-1)^{|y_2||y_1|}[\beta(y_2), [y_3, y_1]_\Gamma]_\Gamma)\\
&=& (0, 0).
\end{eqnarray*}
Also, $(L\oplus \Gamma, [~,~]_{L\oplus \Gamma}, \gamma)$ is multiplicative, since for all homogeneous elements $(x_1, y_1), (x_2, y_2)\in L\oplus \Gamma$, we have
\begin{eqnarray*}
\gamma([(x_1, y_1), (x_2, y_2)]_{L\oplus \Gamma})&=&\gamma([x_1, x_2]_L, [y_1, y_2]_{\Gamma})\\
                                                 &=& ([\alpha(x_1), \alpha(x_2)]_L, [\beta(y_1), \beta(y_2)]_{\Gamma})\\
                                                 &=&([(\alpha(x_1), \beta(y_1)), (\alpha(x_2), \beta(y_2))]_{L\oplus \Gamma}\\
                                                 &=&[\gamma(x_1, y_1), \gamma(x_2, y_2)]_{L\oplus \Gamma}.
\end{eqnarray*}
Finally,
\begin{eqnarray*}
(L\oplus \Gamma)^0&=&\{(x,y)\in(L\oplus \Gamma)_0~|~\gamma(x, y)=(x,y)\}\\
&=&\{(x, y)\in (L\oplus \Gamma)_0~|~(\alpha(x), \beta(y))=(x, y)\}\\
&=&\{(x, y)\in(L\oplus \Gamma)_0~|~\alpha(x)=x~\mathrm{and}~\beta(y)=y\}\\
&=& L^0\oplus \Gamma^0.
\end{eqnarray*}
\hfill $\Box$
\begin{theorem}\label{A9}
Given two restricted Hom-Lie superalgebras $(L, [~, ~]_L, \alpha, [p]_1)$ and $(\Gamma, [~,~]_{\Gamma}, \beta, [p]_2)$, then $(L\oplus \Gamma, [~,~]_{L\oplus \Gamma}, \gamma, [p])$ has a restricted Hom-Lie superalgebras, where $[~,~]_{L\oplus \Gamma}$ and $\gamma$ are defined as in Theorem \ref{A8}, and
$$[p]: (L\oplus \Gamma)^0\rightarrow (L\oplus \Gamma)^0;~(u,v)\mapsto (u^{[p]_1}, v^{[p]_2}).$$
\end{theorem}
{\it Proof.~}
According to Theorem \ref{A8}, it is enough to check the three conditions given in Definition \ref{A2}. Let $(x_1, y_1)\in L^0\oplus \Gamma^0$ and $(x_2, y_2)\in L\oplus \Gamma$. We have
\begin{eqnarray*}
(\mathrm{ad}_{\gamma}(x_1, y_1))^p(x_2, y_2)&=& (\mathrm{ad}_{\gamma}(x_1, y_1))^{p-1}(\mathrm{ad}_\gamma (x_1, y_1))(x_2, y_2)\\
&=&(\mathrm{ad}_\gamma (x_1, y_1))^{p-1}[(x_1, y_1), (\alpha(x_2), \beta(y_2))]_{L\oplus \Gamma}\\
&=&(\mathrm{ad}_\gamma(x_1, y_1))^{p-1}([x_1, \alpha(x_2)]_L, [y_1, \beta(y_2)]_{\Gamma})\\
&=&(\mathrm{ad}_\gamma(x_1, y_1))^{p-2}\left([x_1, [x_1, \alpha^2(x_2)]_L]_L, [y_1, [y_1, \beta^2(y_2)]_{\Gamma}]_{\Gamma}\right)\\
&\vdots&\\
&=&\left([x_1, [x_1, \cdots, [x_1, \alpha^p(x_2)]_L\cdots]_L]_L, [y_1, [y_1, \cdots, [y_1, \beta^p(y_2)]_{\Gamma}\cdots]_{\Gamma}]_{\Gamma}\right)\\
&=&\left(\left(\mathrm{ad}_\alpha x_1\right)^p(x_2), \left(\mathrm{ad}_{\beta} y_1\right)^p(y_2)\right)\\
&=&\left(\mathrm{ad}_{\alpha}x_1^{[p]_1}(x_2), \mathrm{ad}_{\beta}y_1^{[p]_2}(y_2)\right)\\
&=&\left([x_1^{[p]_1}, \alpha(x_2)]_L, [y_1^{[p]_2}, \beta(y_2)]_{\Gamma}\right)\\
&=&\left[\left(x_1^{[p]_1}, y_1^{[p]_2}\right), \left(\alpha(x_2), \beta(y_2)\right)\right]_{L\oplus \Gamma}\\
&=&\mathrm{ad}_{\gamma}(x_1, y_1)^{[p]}(x_2, y_2).
\end{eqnarray*}
This proves that $(\mathrm{ad}_{\gamma}(x_1, y_1))^{[p]}=\mathrm{ad}_{\gamma}(x_1, y_1)^{[p]}$. For $k\in F, (x, y)\in L^0\oplus \Gamma^0$, we obtain
$$\left(k(x, y)\right)^{[p]}=\left((kx)^{[p]_1}, (ky)^{[p]_2}\right)=\left(k^px^{[p]_1}, ky^{[p]_2}\right)=k^p\left(x^{[p]_1}, y^{[p]_2}\right)= k^p(x, y)^{[p]}.$$
Finally, for $(x_1, y_1), (x_2, y_2)\in L^0\oplus \Gamma^0$, one gets
\begin{eqnarray*}
\left((x_1, y_1)+(x_2, y_2)\right)^{[p]}&=& \left(x_1+x_2, y_1+y_2\right)^{[p]}\\
&=&\left((x_1+x_2)^{[p]_1}, (y_1+y_2)^{[p]_2}\right)\\
&=&\left(x_1^{[p]_1}+x_2^{[p]_1}+\sum_{i=1}^{p-1}s_i(x_1, x_2), y_1^{[p]_2}+y_2^{[p]_2}+\sum_{i=1}^{p-1}s_i(y_1, y_2)\right)\\
&=&\left(x_1^{[p]_1}, y_1^{[p]_2}\right)+\left(x_2^{[p]_1}, y_2^{[p]_2}\right)+ \left(\sum_{i=1}^{p-1}s_i(x_1, x_2), \sum_{i=1}^{p-1}s_i(y_1, y_2)\right)\\
&=& (x_1, y_1)^{[p]}+ (x_2, y_2)^{[p]}+\sum_{i=1}^{p-1}s_i\left((x_1, y_1), (x_2, y_2)\right).
\end{eqnarray*}
\hfill $\Box$
\begin{corollary}\label{A10}
Suppose that $\left(L_1, [~,~]_{L_1}, \alpha_1. [p]_1\right), \ldots, \left(L_n, [~,~]_{L_n}, \alpha_n, [p]_n\right)$ are restricted Hom-Lie superalgebras. Then $\left(L_1\oplus\cdots\oplus L_n, [~,~]_{L_1\oplus\cdots\oplus L_n}, \gamma, [p]\right)$ is a restricted Hom-Lie superalgebra, where the bilinear map $[~,~]_{L_1\oplus\cdots\oplus L_n}$ is defined on homogenous elements by
$$\left[(x_1, \ldots, x_n), (y_1, \ldots, y_n)\right]_{L_1\oplus\cdots\oplus L_n}=\left([x_1, y_1]_{L_1}, \ldots, [x_n, y_n]_{L_n}\right),$$
and the linear map $\gamma:L_1\oplus\cdots\oplus L_n\rightarrow L_1\oplus\cdots\oplus L_n$ is defined as
$$\gamma(x_1, \ldots, x_n)=(\alpha_1(x_1), \ldots, \alpha_n(x_n)).$$
Also, the $p$-map $[p]: (L_1\oplus \cdots\oplus L_n)^0\rightarrow (L_1\oplus\cdots\oplus L_n)^0$ is given by
$$\left(x_1, \ldots, x_n\right)^{[p]}=\left(x_1^{[p]_1}, \ldots, x_n^{[p]_n}\right).$$
\end{corollary}
\section{On morphisms of Hom-Lie superalgebras}

\begin{definition}[{{\cite{Makhlouf}}}]\label{A6}
Let $(L, [~,~]_L, \alpha)$ and $(\Gamma, [~,~]_\Gamma, \beta)$ be two Hom-Lie superalgebras. An even superspace $\varphi: L \rightarrow \Gamma$ is said to be a morphism of Hom-Lie superalgebras if
\begin{eqnarray*}
[\varphi(u), \varphi(v)]_\Gamma=\varphi([u, v]_L)~\forall u, v\in L~\mathrm{and}~\varphi\circ \alpha=\beta\circ \varphi.
\end{eqnarray*}
\end{definition}
If $\varphi: L\rightarrow \Gamma$ is an even superspace, then the graph of $\varphi$ is the set
$$G_\varphi = \{(x, \varphi(x)~|~x\in L\}\subseteq L\oplus \Gamma.$$
\begin{definition}\label{A7}
A morphism of Hom-Lie superalgebra
$$\varphi: (L, [~,~]_L, \alpha, [p]_1)\rightarrow (\Gamma, [~,~]_\Gamma, \beta, [p]_2)$$
is said to be restricted if $\varphi\left(x^{[p]_1}\right)=(\varphi(x))^{[p]_2}~\forall x\in L^0$.
\end{definition}
\begin{theorem}\label{A11}
A linear map $\varphi: (L, [~,~]_L, \alpha, [p]_1)\rightarrow (\Gamma, [~,~]_{\Gamma}, \beta, [p]_2)$ is a restricted morphism of Hom-Lie superalgebras if and only if the graph $G_{\varphi}$ is a restricted Hom-Lie subsuperalgebra of $(L\oplus\Gamma, [~,~]_{L\oplus\Gamma}, \gamma, [p])$ where $[~,~]_{L\oplus \Gamma}, \gamma$, and $[p]$ are defined as in Theorem \ref{A9}.
\end{theorem}
{\it Proof.~}
Suppose $\varphi:(L, [~,~]_L, \alpha, [p]_1)\rightarrow (\Gamma, [~,~]_{\Gamma}, \beta, [p]_2)$ is a restricted morphism of Hom-Lie superalgebras. Let $(x_1, \varphi(x_1)), (x_2, \varphi(x_2))\in G_{\varphi}$. Then
\begin{eqnarray*}
[(x_1, \varphi(x_1)), (x_2, \varphi(x_2))]_{L\oplus \Gamma}&=& \left([x_1, x_2]_L, [\varphi(x_1), \varphi(x_2)]_{\Gamma}\right)\\
&=& \left([x_1, x_2]_L, \varphi([x_1, x_2])\right)\in G_{\varphi}.
\end{eqnarray*}
This shows $G_\varphi$ is closed under the bracket operation. For $(x_1, \varphi(x_1))\in G_\varphi$, we have
$$\gamma(x_1, \varphi(x_1))=(\alpha(x_1), \beta(\varphi(x_1)))= (\alpha(x_1), \varphi(\alpha(x_1)))\in G_\varphi.$$
So far we have shown that $G_{\varphi}$ is a Hom-Lie subsuperalgebra. For $(x_1, \varphi(x_1))\in G_{\varphi}^0$, we have
\begin{eqnarray*}
\left(x_1, \varphi(x_1)\right)^{[p]}&=&\left(x_1^{[p]_1}, \left(\varphi(x_1)\right)^{[p]_2}\right)\\
&=&\left(x_1^{[p]_1}, \varphi\left(x_1^{[p]_1}\right)\right)\in G_{\varphi}^0.
\end{eqnarray*}
Conversely, suppose $G_{\varphi}$ is a restricted Hom-Lie subsuperalgebra of $(L\oplus \Gamma, [~,~]_{L\oplus \Gamma}, \gamma, [p])$. Let $x_1, x_2\in L$. Then $\left[(x_1, \varphi(x_1)), (x_2, \varphi(x_2))\right]_{L\oplus \Gamma}=\left([x_1, x_2]_L, [\varphi(x_1), \varphi(x_2)]_{\Gamma}\right)\in G_{\varphi}$, and so $\varphi\left([x_1, x_2]_L\right)= \left[\varphi(x_1), \varphi(x_2)\right]_{\Gamma}$. For $x\in L$, we have $\gamma(x, \varphi(x))= (\alpha(x), \beta(\varphi(x)))\in G_{\varphi}$, so that $\varphi(\alpha(x))= \beta(\varphi(x))$. Finally, if $x\in L^0$, then we have $(x, \varphi(x))^{[p]}=(x^{[p]_1}, (\varphi(x))^{[p]_2})\in G_{\varphi}$, and so $\varphi(x^{[p]_1})= (\varphi(x))^{[p]_2}$. \hfill $\Box$

We extend, in the following theorem, the result in \cite{Makhlouf} to restricted Lie superalgebras case. It gives a way to construct restricted Hom-Lie superalgebras starting from a restricted Lie superalgebra and an even superalgebra endomorphism.
\begin{theorem}\label{AA11}
Let $(L, [~,~], [p])$ be a restricted Lie superalgebra and $\alpha: L\rightarrow L$ be an even Lie superalgebra endomorphism. Then $(L, [~,~]_{\alpha}, \alpha, [p]|_{L^0})$, where $[x, y]_{\alpha}=\alpha([x, y])$ is a restricted Hom-Lie superalgebra.
\end{theorem}
{\it Proof.~} The bracket is obviously super skew-symmetric. Let $x, y, z\in L$ be homogenous elements. With a direct computation, we have
\begin{eqnarray*}
&&(-1)^{|x||z|}[\alpha(x), [y,z]_{\alpha}]_{\alpha}+(-1)^{|z||y|}[\alpha(z), [x, y]_{\alpha}]_{\alpha}+(-1)^{|y||x|}[\alpha(y), [z,x]_{\alpha}]_{\alpha}\\
&=& (-1)^{|x||z|}\alpha^2\left([x, [y,z]]\right)+(-1)^{|z||y|}\alpha^2\left([z, [x, y]]\right)+(-1)^{|y||x|}\alpha^2\left([y, [z,x]]\right)\\
&=&\alpha^2\left((-1)^{|x||z|}[x, [y,z]]+(-1)^{|z||y|}[z, [x, y]]+(-1)^{|y||x|}[y, [z,x]]\right)=0.
\end{eqnarray*}
Now, the Hom-Lie superalgebra $(L, [~,~]_{\alpha}, \alpha)$ is multiplicative. Indeed, for homogeneous elements $x, y\in L$, we have $\alpha\left([x, y]_{\alpha}\right)= \alpha^2([x, y])=[\alpha(x), \alpha(y)]_{\alpha}$.\\
Next, we check the three conditions given in Definition \ref{A2}. For $x\in L^0$ and $z\in L$, we have
\begin{eqnarray*}
\left(\mathrm{ad}_{\alpha}x\right)^p(z)&=&\left(\mathrm{ad}_{\alpha}x\right)^{p-1}\left(\mathrm{ad}_{\alpha}x\right)(z)\\
&=&\left(\mathrm{ad}_{\alpha}x\right)^{p-1}\left([x, \alpha(z)]_{\alpha}\right)\\
&=&\left(\mathrm{ad}_{\alpha}x\right)^{p-1}([x, \alpha^2(z)])\\
&=&\left(\mathrm{ad}_{\alpha}x\right)^{p-2}\left([x, [x, \alpha^2(z)]]_{\alpha}\right)\\
&=&\left(\mathrm{ad}_{\alpha}x\right)^{p-2}([x, [x, \alpha^3(z)]])\\
&\vdots&\\
&=&[x, [x, \cdots, [x, \alpha^{p+1}(z)]\cdots]]\\
&=&\left(\mathrm{ad}x\right)^p(\alpha(z))\\
&=&\mathrm{ad}\left(x^{[p]}\right)(\alpha(z))\\
&=&\left[x^{[p]}, \alpha(z)\right]\\
&=&\mathrm{ad}_{\alpha}x^{[p]|_{L^0}}(z).
\end{eqnarray*}
The second and the third properties are clear.\hfill $\Box$
\begin{theorem}\label{AA12}
Let $(L, [~,~]_L, [p]_1)$ and $(\Gamma, [~,~]_{\Gamma}, [p]_2)$ be restricted Lie superalgebras, $\alpha:L\rightarrow L$ and $\beta:\Gamma\rightarrow \Gamma$ be an even Lie superalgebra endomorphisms, and $f:L\rightarrow \Gamma$ be a morphism of restricted Lie superalgebras with $f\circ \alpha= \beta\circ f$. Then
$$f: (L, [~,~]_{\alpha}, \alpha, [p]_1)\rightarrow (\Gamma, [~,~]_{\beta}, \beta, [p]_2),$$
where  $[~,~]_{\alpha}$ and $[~,~]_{\beta}$ are defined as in Theorem \ref{AA11}, is a morphism of restricted Hom-Lie superalgebras.
\end{theorem}
{\it Proof.~}
If $u, v\in L$, we have
$$f\left([u, v]_{\alpha}\right)= f(\alpha[u, v])=\beta\circ f([u, v])=\beta([f(u), f(v)])=[f(u), f(v)]_{\beta}.$$
Also. for $u\in L^0$ and using $L^0\subseteq L_0$, we obtain $f\left(u^{[p]_1}\right)= \left(f(u)\right)^{[p]_2}$.\hfill $\Box$
\begin{theorem}\label{AA13}
Let $(L, [~,~]_L, \alpha)$ and $(\Gamma,[~,~]_{\Gamma}, \beta)$ be multiplicative Hom-Lie superalgebras, and $\varphi$ be a one to one morphism of Hom-Lie superalgebras. If $C$ is a Hom-Lie subsuperalgebra of $\Gamma$ and $[p]: C^0\rightarrow C^0$ is a $p$-map, then $\left(\varphi^{-1}(C), [~,~]_L, \alpha|_{\varphi^{-1}(C)}, [p]'\right)$, where
$$[p]': \left(\varphi^{-1}(C)\right)^0\rightarrow \left(\varphi^{-1}(C)\right)^0; ~x\mapsto \varphi^{-1}\left(\left(\varphi(x)\right)^{[p]}\right),$$
is a restricted Hom-Lie superalgebra.
\end{theorem}
{\it Proof.~}
First, we show $\varphi^{-1}(C)$ is a Hom-Lie subsuperalgebra of $L$. For $x_1, x_2\in \varphi^{-1}(C)$, there exist $y_1, y_2\in C$ with $\varphi(x_1)=y_1$ and $\varphi(x_2)=y_2$. Now, $\varphi([x_1, x_2])=[y_1, y_2]$, implies $[x_1, x_2]\in \varphi^{-1}(C)$. Also, $\alpha(x_1)\in \varphi^{-1}(C)$ follows from $\beta(\varphi(x_1))=\varphi(\alpha(x_1))$. Next, if $x\in \left(\varphi^{-1}(C)\right)^0$, then it is clear that $\varphi(x)\in C^0$. Also, for $k\in F$, we have
\begin{eqnarray*}
(kx)^{[p]'}&=& \varphi^{-1}\left((\varphi(kx))^{[p]}\right)\\
&=&\varphi^{-1}\left((k\varphi(x))^{[p]}\right)\\
&=&\varphi^{-1}\left(k(\varphi(x))^{[p]}\right)\\
&=& k\left(\varphi(x)\right)^{[p]}= kx^{[p]'}.
\end{eqnarray*}
Also, for $x_1, x_2\in \left(\varphi^{-1}(C)\right)^0$ we have
\begin{eqnarray*}
(x_1+x_2)^{[p]'}&=& \varphi^{-1}\left(\varphi(x_1)+\varphi(x_2)\right)^{[p]}\\
&=& \varphi^{-1}\left(\varphi(x_1)+\varphi(x_2)\right)^{[p]}\\
&=& \varphi^{-1}\left((\varphi(x_1))^{[p]}+(\varphi(x_2))^{[p]}+\sum_{i=1}^{p-1}s_i(\varphi(x_1), \varphi(x_2))\right)\\
&=&\varphi^{-1}\left((\varphi(x_1))^{[p]}\right)+ \varphi^{-1}\left((\varphi(x_2))^{[p]}\right) + \varphi^{-1}\left(\sum_{i=1}^{p-1}s_i(\varphi(x_1), \varphi(x_2))\right)\\
&=& x_1^{[p]'}+ x_2^{[p]'}+ \sum_{i=1}^{p-1}s_i(x_1, x_2).
\end{eqnarray*}
Finally, let $x_1\in \left(\varphi^{-1}(C)\right)^0$ and $x_2\in \phi^{-1}(C)$. There exist $y_1\in C^0$ and $y_2\in C$ with $\varphi(x_1)= y_1$ and $\varphi(x_2)=y_2$. With a direct computation, we have
\begin{eqnarray*}
\mathrm{ad}_{\alpha}x^{[p]'}(x_2)&=& \left[x_1^{[p]'}, \alpha(x_2)\right]\\
&=&\left[\varphi^{-1}\left(\left(\varphi(x_1)\right)^{[p]}\right), \alpha(x_2)\right]\\
&=&\varphi^{-1}\left[\left(\varphi(x_1)\right)^{[p]}, \varphi\circ \alpha(x_2)\right]\\
&=& \varphi^{-1}\left[y_1^{[p]}, \beta(y_2)\right]\\
&=& \varphi^{-1}\left(\mathrm{ad}_{\beta}y_1^{[p]}\right)(y_2)\\
&=&\varphi^{-1}\left(\mathrm{ad}_{\beta}y_1\right)^{p}(y_2)\\
&=&\varphi^{-1}\left(\mathrm{ad}_{\beta}y_1\right)^{p-1}[\varphi(x_1), \beta(\varphi(x_2))]\\
&=&\varphi^{-1}\left(\mathrm{ad}_{\beta}y_1\right)^{p-1}[\varphi(x_1), \varphi(\alpha(x_2))]\\
&\vdots&\\
&=&\varphi^{-1}[\varphi(x_1), [\varphi(x_1), [\cdots, [\varphi(x_1), \varphi(\alpha^p(x_2))]\cdots]]]\\
&=&[x_1, [x_1, \cdots, [x_1, \alpha^p(x_2)]\cdots]]\\
&=&\left(\mathrm{ad}_{\alpha}x_1\right)^p(x_2).
\end{eqnarray*}
\hfill $\Box$
\begin{theorem}\label{AA14}
Suppose that $(L, [~,~]_L, \alpha, [p])$ is a restricted Hom-Lie superalgebra, $(\Gamma, [~,~]_{\Gamma}, \beta)$ is a multiplicative Hom-Lie superalgebra, and $\varphi: L\rightarrow \Gamma$ is a one to one morphism of Hom-Lie superalgebras. Then $(\varphi(L), [~,~]_{\Gamma}, \beta|_{\varphi(L)}, [p]')$ is a restricted Hom-Lie superalgebra, where
$$[p]': (\varphi(L))^0\rightarrow (\varphi(L))^0;~\varphi(x)\mapsto \varphi\left(x^{[p]}\right).$$
\end{theorem}
{\it Proof.~} For $\varphi(x_1), \varphi(x_2)\in \varphi(L)$, we have $[\varphi(x_1), \varphi(x_2)]= \varphi([x_1, x_2])\in \varphi(L)$ and $\beta(\varphi(x_1))=\varphi(\alpha(x_1))\in \varphi(L)$. For $\varphi(x)\in (\varphi(L))^0$, we have $\beta(\varphi(x))= \varphi(\alpha(x))=\varphi(x)$. Since $\varphi$ is one to one we have $\alpha(x)= x$. Thus, $x\in L^0$. Next, we check the three conditions given in Definition \ref{A2}. For $\varphi(x)\in (\varphi(L))^0$ and $k\in F$, we obtain
$$(k\varphi(x))^{[p]'}= (\varphi(kx))^{[p]'}=\varphi\left((kx)^{[p]}\right)=\varphi\left(kx^{[p]}\right)=k\varphi\left(x^{[p]}\right)= k(\varphi(x))^{[p]'}.$$
If $\varphi(x_1), \varphi(x_2)\in (\varphi(L))^0$, then we have
\begin{eqnarray*}
(\varphi(x_1)+ \varphi(x_2))^{[p]'}&=& (\varphi(x_1+x_2))^{[p]'}\\
&=&\varphi\left((x_1+x_2)^{[p]}\right)\\
&=& \varphi\left(x_1^{[p]}+ x_2^{[p]}+\sum_{i=1}^{p-1}s_i(x_1, x_2)\right)\\
&=&\varphi\left(x_1^{[p]}\right)+ \varphi\left(x_2^{[p]}\right)+ \varphi\left(\sum_{i=1}^{p-1}s_i(x_1, x_2)\right)\\
&=&\left(\varphi(x_1)\right)^{[p]'}+ \left(\varphi(x_2)\right)^{[p]'}+ \sum_{i=1}^{p-1}s_i(\varphi(x_1), \varphi(x_2)).
\end{eqnarray*}
In addition, for $\varphi(x_1)\in (\varphi(L))^0$ and $\varphi(x_2)\in \varphi(L)$, we obtain
\begin{eqnarray*}
\mathrm{ad}_{\beta}(\varphi(x_1))^{[p]'}(\varphi(x_2))&=& \left[(\varphi(x_1))^{[p]'}, \beta(\varphi(x_2))\right]\\
&=&\left[(\varphi(x_1))^{[p]'}, \varphi(\alpha(x_2))\right]\\
&=&\varphi\left(\left[x_1^{[p]}, \alpha(x_2)\right]\right)\\
&=& \varphi\left(\mathrm{ad}_{\alpha}x_1^{[p]}\right)(x_2)\\
&=& \varphi\left(\left(\mathrm{ad}_{\alpha}x_1\right)^p(x_2)\right)\\
&=& \varphi\left([x_1, [x_1, \cdots, [x_1, \alpha^p(x_2)]\cdots]]\right)\\
&=&[\varphi(x_1), [\varphi(x_1), \cdots, [\varphi(x_1), \varphi(\alpha^p(x_2))]\cdots]]\\
&=&[\varphi(x_1), [\varphi(x_1), \cdots, [\varphi(x_1), \beta^p(\varphi(x_2))]\cdots]]\\
&=&\left(\mathrm{ad}_{\beta}\varphi(x_1)\right)^p(\varphi(x_2)).
\end{eqnarray*}
\hfill $\Box$
\section{Restrictable Hom-Lie superalgebra}
\begin{definition}\label{A20}
A multiplicative Hom-Lie superalgebra $(L, [~,~], \alpha)$ is called restrictable if $(\mathrm{ad}_{\alpha}x)^p\in \mathrm{ad}_{\alpha}L^0$ for all $x\in L^0$, where $\mathrm{ad}_{\alpha}L^0=\{\mathrm{ad}_{\alpha}x~|~x\in L^0\}$.
\end{definition}

The following theorem was obtained by Guan and Chen in \cite{Guan} in the setting of Lie algebras. We extend it to Lie superalgebra case.

\begin{theorem}\label{A21}
A multiplicative Hom-Lie superalgebra is restrictable if and only if there is a $p$-map $[p]: L^0\rightarrow L^0$ which makes $L$ a restricted Hom-Lie superalgebra.
\end{theorem}
{\it Proof.~}
Suppose $(L, [~,~], \alpha, [p])$ is a restricted Hom-Lie superalgebra. Let $x\in L^0$. Then $(\mathrm{ad}_{\alpha}x)^p=\mathrm{ad}_{\alpha}x^{[p]}\in \mathrm{ad}_{\alpha}L^0$. Conversely, suppose that $(L, [~,~], \alpha)$ is restrictable. Let $\{e_j~|~j\in J\}$ be a basis of $L^0$. Then for each $j\in J$, there exists $y_j\in L^0$ such that $(\mathrm{ad}_{\alpha}e_j)^p=\mathrm{ad}_{\alpha}y_j$. Now, $e_j^p-y_j\in C(L^0)$, where $$e_j^p=[e_j, [e_j,\cdots,[e_j, e_j]\cdots]],$$
for if $z\in L^0$, then $(\mathrm{ad}_{\alpha}e_j)^p(z)-\mathrm{ad}_{\alpha}y_j(z)= [e_j^p-y_j, \alpha(z)]=0$. Define a function
$$f: L^0\rightarrow C(L^0);~ \sum c_ie_i\mapsto \sum c_i^p(y_i-e_i^p).$$
Clearly $f$ is a $p$-semilinear map. Set
$$W=\{x\in L^0~|~x^p+f(x)\in L^0\}.$$
Then $W$ is a subspace of $L^0$. Indeed, if $x, y\in W$ and $a, b\in F$, then
$$(ax+by)^p+f(ax+by)=a^px^p+b^py^p+\sum_{i=1}^{p-1}s_i(ax, by)+k^pf(x)+f(by)\in L^0.$$
Since $e_j^p+f(e_j)\in L^0$, it follows $x^p+ f(x)\in L^0$ for all $x\in L^0$. By Theorem \ref{A4}, we have
$$[p]: L^0\rightarrow L^0;~x\mapsto x^{p}+f(x)$$
is a $p$-map with $e_j^{[p]}=y_j$.\hfill $\Box$

Let $(L, [~,~], \alpha)$ be a Hom-Lie superalgebra, and let $U$ and $W$ be subspaces of $L$. Then $L$ is a direct sum of $U$ and $W$, and we  write $L=U\oplus W$, if $L=U+W$ and $U\cap W=\{0\}$. The subspace $U$ is a Hom-Lie ideal of $L$ if $\alpha(U)\subseteq U$ and $[x, y]\in U$ for all $x\in U$ and $y\in L$. We have the following result.
\begin{theorem}\label{A23}
Let $(L, [~,~], \alpha)$ be a Hom-Lie superalgebra and let $U$ and $W$ be Hom-Lie ideals of $L$ with $L=U\oplus W$. Then, $L$ is restrictable if and only if $U$ and $W$ are restrictable.
\end{theorem}
{\it Proof.~}
Suppose $L$ is restrictable. Then by Theorem \ref{A21}, we have $L$ is restricted. Now, the result follows from the trivial fact that a subalgebra of a restricted Hom-Lie superalgebra is restricted and Theorem \ref{A21}. Conversely, suppose that $U$ and $V$ are restrictable. Let $x=x_1+ x_2$, where $x_1\in U$ and $x_2\in W$, be an element of $L^0$. Then $\alpha(x)=\alpha(x_1)+\alpha(x_2)=x_1+ x_2$, and so $x_1\in U^0$ and $x_2\in W^0$. As $U$ and $W$ are restrictable, there exist $y_1\in Y^0$ and $y_2\in W^0$ such that $(\mathrm{ad}_{\alpha}x_1)^p=\mathrm{ad}_{\alpha}y_1$ and $(\mathrm{ad}_{\alpha}x_2)^p=y_2$. Now,
\begin{eqnarray*}
(\mathrm{ad}_{\alpha}(x_1+x_2))^p&=&(\mathrm{ad}_{\alpha}x_1+\mathrm{ad}_{\alpha}x_2)^p\\
&=&(\mathrm{ad}_{\alpha}x_1)^p+(\mathrm{ad}_{\alpha}x_2)^p\\
&=&\mathrm{ad}_{\alpha} y_1+\mathrm{ad}_{\alpha}y_2\\
&=&\mathrm{ad}_{\alpha}(y_1+y_2).
\end{eqnarray*}
Therefore, $L$ is restrictable.\hfill $\Box$


\end{document}